



\magnification=\magstep 1
\parindent = 0 pt
\baselineskip = 16 pt
\parskip = \the\baselineskip

\font\AMSBoldBlackboard = msbm10

\def\QQ{{\hbox{\AMSBoldBlackboard Q}}}

\def\ZZ{{\hbox{\AMSBoldBlackboard Z}}}

\settabs 12\columns
\rightline{math.NT/9901051}
\rightline{v1 12-JAN-99}
\rightline{v2 31-JAN-99}
\vskip 1 true in
{\bf \centerline{SCATTERING ON THE P-ADIC FIELD AND A TRACE FORMULA}}
\vskip 0.5 true in
\centerline{Jean-Fran\c{c}ois Burnol}
\par
\centerline{January 1999}
\par
I apply the set-up of Lax-Phillips Scattering Theory to a non-archimedean local field. It is possible to choose the outgoing space and the incoming space to be Fourier transforms of each other. Key elements of the Lax-Phillips theory are seen to make sense and to have the expected interrelations: the scattering matrix S, the projection K to the interacting space, the contraction semi-group Z and the time delay operator T. The scattering matrix is causal, its analytic continuation has the expected poles and zeros, and its phase derivative is the (non-negative) spectral function of T, which is also the restriction to the diagonal of the kernel of K. The contraction semi-group Z is related to S (and T) through a trace formula. Introducing an odd-even grading on the interacting space allows to express the Weil local explicit formula in terms of a ``supertrace''. I also apply my methods to the evaluation of a trace considered by Connes.
\vfill
{\parskip = 0 pt
62 rue Albert Joly\par
F-78000 Versailles\par
France\par}
\eject

{\bf
TABLE OF CONTENTS\par
\par
Introduction\par
Additive and multiplicative Fourier analysis on a local field\par
The scattering matrix\par
The interacting space and its associated kernel\par
The time delay operator and the conductor operator\par
The contraction semi-group and a trace formula\par
Weil's local term in terms of a supertrace\par
Exact evaluation of a trace considered by Connes\par
References\par}
\par
\vfill \eject

{\bf Introduction}

It is well-known (Wiener) that the Fourier Transform on the real line is obtained at a phase of $\pi/2$ in the quantum time evolution of the harmonic oscillator (with suitably chosen parameters, and with the zero-point energy removed). Unfortunately it seems difficult to make sense in a really satisfactory way of the infinitesimal generator of this group on a $p$-adic local field.

It is less well-known that the composition of the Fourier Transform with the Inversion ($\varphi(x) \mapsto {1\over|x|}\varphi({1\over x})$), which is a dilation invariant operator is also obtained as a member of a unitary representation, this time non-compact, of the additive group of the reals. This follows from the existence of the logarithm of the Tate-Gel'fand-Graev Gamma function on the critical line. Now this makes sense not only for the reals, but also for the complex numbers, and in fact also for any completion of an algebraic number field at a finite non-ramified place. One has to restrict oneself in that case to that part of the Hilbert space invariant under rotations by units, and an explicit computation shows that the Tate-Gel'fand-Graev Gamma function on the critical line has a logarithm respecting the inherent periodicity in the vertical direction (the next chapter contains all the necessary set-up to understand this sentence). 

But is there a ``physical'' counterpart to these observations as there is one, indeed very basic to quantum physics, for the Wiener Theorem about fractional powers of the Fourier Transform? It seems that the question has much more to do with scattering processes than with bound states problems, and I was thus led to try to apply one of the various rigorous mathematical methodologies concerning scattering to p-adic numbers.

I will show in this paper that it is possible to a surprising extent to adapt the axiomatic set-up of Lax-Phillips [LP90] to a very basic situation involving the non-archimedean Fourier Transform, and to translate almost everything from their paper [LP78] to the non-archimedean field. Of course the analysis is much easier, as we will compute only one special, but very interesting case, and will compute everything completely explicitely. In the end we obtain a (super-)trace formula giving Weil's local term of the Explicit Formula of analytic number theory. I have no deeper reason to offer for the introduction of the grading necessary to get Weil's local term, than the fact that it does the job. 

Another motivation is to establish a bridge between some computations and results of Connes [Co98] and my own considerations on the conductor operator and the Explicit Formula [Bu98a,b,c]. In the last section I evaluate exactly a trace considered by Connes (Theorem $3$ of Section V of [Co98]).

{\bf Additive and multiplicative Fourier analysis on a local field}

Let us start with some notations concerning the Fourier Analysis on local fields [Ta50] and previous results of the author [Bu98a,b,c].

Let $K$ be a number field and $\nu$ one of its places with completion $K_\nu$. In this paper only the case of a non-archimedean place will be considered. Let $p$ be the prime number over which lies $\nu$ and $q$ the cardinality of the residue field (which is a power of $p$). We will write $x$ or $y$ for an element of $K_\nu$, viewed additively, and $t$ or $u$ for an element of the multiplicative group $K_\nu^\times$. To define Fourier transforms we choose the basic additive character to be $\lambda(x) = \exp(2\pi\ i\ Tr(x))$, where ``Tr'' is the trace down to $\QQ_p$. The Fourier Transform ${\cal F}(\varphi)$ of $\varphi$ (also denoted $\widetilde{\varphi}$) is defined through ${\cal F}(\varphi)(x) = \int_{K_\nu}\varphi(y)\lambda(-xy)\, dy$ where $dy$ is the unique self-dual additive Haar measure corresponding to the choice of $\lambda$ (so that ${\cal F}{\cal F}(\varphi)(y) = \varphi(-y)$). The Fourier Transform is a unitary of the Hilbert space $L^2$ of square-integrable functions with respect to $dy$. We use the bra-ket notation: $<\psi|\varphi> = \int_{K_\nu}\overline{\psi(x)}\varphi(x)\,dx$. The volume of the subring of integers $O_\nu$ is then $q^{-\delta/2}$ where $q$ is the cardinality of the residue field and $\delta$ a non-negative integer called the differental exponent. We choose a uniformizer $\pi$ of $O_\nu$ and normalize the norm according to $|\pi|_\nu = {1\over q}$. One then has $d(ty) = |t|dy$ (we will drop most $\nu$ subscripts as we are dealing only with one place at a time). The characteristic function $\omega$ of $O$ has Fourier transform $\widetilde{\omega}(x) = q^{-\delta/2}\omega(\pi^\delta x)$. There is a unique multiplicative Haar measure on $K_\nu^\times$ which assigns to the compact open subgroup of units a volume of $1$, we denote it by $d^*t$. It is related to the additive Haar measure $dx$ by the formula $a^2 |x|d^*x = dx$ with $a = \sqrt{q^{-\delta/2}(1 - {1\over q})}$. We also use the multiplicative Haar measure $d^\times t = \log(q) d^*t$ which assigns a volume of $\log(q)$ to the units. 

Let's now review briefly the technique used in [Bu98b,c], which will be very useful here.
First the Hilbert space $L^2$ can also be realized as the Hilbert space of square integrable functions on the multiplicative field with respect to $d^*t$ (one assigns to a function $\varphi(x)$ on $K_\nu$ the function $f(t) = a|t|^{1/2}\varphi(t)$ on $K_\nu^\times$). The locally compact abelian group $K_\nu^\times$ has a (unitary) dual $X_\nu$ which is a (countable) union of circles indexed by the characters of the multiplicative group of units. Let $\chi$ be one such character, we still denote by the same letter the unique extension to $K_\nu^\times$ taking the value $1$ at $\pi$. We write $e(\chi)$ for the conductor exponent of $\chi$, which is a non-negative integer vanishing only for the trivial character. The elements of $X_\nu$ in the same component as $\chi$ are then all of the form $t \mapsto \chi(t)|t|^{i\tau}$ for some (real number) $\tau$, and parametrized uniquely by $z = q^{i\tau}$ (the value taken at $1\over\pi$). We note $X_\chi$ for this copy of $U(1)$ and $(\chi; z)$ for the corresponding character of $K_\nu^\times$. The well-known generalization of Fourier Theory applies to this setting and one finds that the measure dual to $d^*t$ reduces on each component to the unique rotation invariant measure of total mass $1$ (we will write ${d\theta\over 2\pi}$ for this measure). We thus have a third incarnation of $L^2$ as $L^2(X,{d\theta\over 2\pi})$, the formula for the isometries being 
$$f(t) \mapsto \widehat{f}(\chi;z) = \int_{K_\nu^\times} f(t)(\chi;z)(t)\,d^*t$$
$$f(t) = \int_X \widehat{f}(\chi;z) \overline{(\chi;z)(t)}\,{d\theta\over 2\pi}$$
Let $U(t)$ be the unitary group of dilations, acting in the additive picture as $\varphi(x) \mapsto |t|^{-1/2}\varphi(x/t)$, in the multiplicative picture as $f(u) \mapsto f(u/t)$, in the spectral picture as $l(\chi;z) \mapsto (\chi;z)(t)\l(\chi;z)$. We will write $L^2_\chi = L^2(X_\chi,{d\theta\over 2\pi})$ for the $\chi$-isotypical component of $L^2$. We note that the action of the units is diagonalized by this decomposition, while the inverse-uniformizer $1\over\pi$ acts as multiplication by the function ``$z$''.

We will now consider various operators. The first example is with the inversion $I$. Defining it to be $\varphi(x) \mapsto {1\over|x|}\varphi({1\over x})$ in the additive picture it becomes $f(u) \mapsto f(1/u)$ in the multiplicative picture and $l(\chi;z) \mapsto l(\overline{\chi};\overline{z})$ in the spectral picture. It is more interesting to look at the additive Fourier Transform ${\cal F}$. The composition ${\cal F}I$ commutes with the dilation action so it should be diagonalized in the spectral picture. Let's call it $\Gamma$. I showed in [Bu98b] how an application of Tate's local functional equation [Ta50] leads to $(\Gamma l)(\chi;z) = \Gamma((\chi;z), {1\over 2})\ l(\chi;z)$ where the $\Gamma$ on the right-hand-side refers to the Tate-Gel'fand-Graev Gamma function on the critical line (almost identical functions are tabulated in [Ta50]). Abbreviating $\Gamma((\chi;z), {1\over 2})$ to $\Gamma(\chi;z)$ this gives explicitely
$$\hbox{for $\chi$ trivial:}\;\ \Gamma(1;z) = z^\delta\ {1-{z\over\sqrt{q}}\over{1-{\overline{z}\over\sqrt{q}}}}$$
$$\hbox{for $\chi$ ramified:}\;\ \Gamma(\chi;z) = w(\chi)\ z^{e(\chi) + \delta}$$
The ``root number'' $w(\chi)$ is a certain complex number of unit modulus. We won't need its value in what follows. For later reference let us recap this as:

{\bf Theorem I:}\quad ${\cal F}(l)(\chi;z) = \Gamma(\chi;z)\ l(\overline{\chi};\overline{z})\;\bullet$

One corollary to this is the existence of an interesting domain in $L^2$, defined in the spectral picture to be the domain of functions smooth on each component and non-vanishing only on finitely many components. This domain, let's call it $\Delta$, is then obviously stable under dilations, inversion, Fourier Transform. Let also ${\cal S}_0(K)$ be the dense domain of Schwartz-Bruhat functions (locally constant with compact support) vanishing at the origin. Finally let ${\cal S}(K)$ be the domain of all Schwartz-Bruhat functions. One has ${\cal S}_0(K)\subset{\cal S}(K)\subset\Delta$.

Let's now consider some unbounded operators. First let $A$ be given in the additive picture as multiplication with $\log(|x|)$. As this is unbounded we should specify a domain $\Delta(A)$, and we take it at first to be ${\cal S}_0(K)$. Obviously it is stable under $A$, and $A+i$ as well as $A-i$ are invertible on it, so $(A,{\cal S}_0(K))$ is essentially self-adjoint. One finds that $A$ acts in the spectral picture on $\Delta$ as the differential operator $D = \log(q) \ z {\partial\over \partial z}$.

The conjugate $B = {\cal F}A{\cal F}^{-1} = {\cal F}^{-1}A{\cal F}$ with domain $\Delta$ is thus essentially self-adjoint. It is easier to evaluate first $IBI = \Gamma^{-1}A\Gamma = A + H$ where $H$ is the dilation invariant operator with spectral function $H(\chi;z) = \Gamma^{-1}(\chi;z)D(\Gamma(\chi;z))$. Explicitely
$$H(1;z) = \delta\log(q) - \left({{z\over\sqrt{q} - z}}+{\overline{z}\over{\sqrt{q} - \overline{z}}}\right)\log(q) = (\delta + 1)\log(q) - \left({1 - {1\over q} \over |1 - {z \over \sqrt{q}}|^2 }\right)\log(q)$$
$$\hbox{for $\chi$ ramified:}\;\ H(\chi;z) = (e(\chi) + \delta)\log(q)$$
From this we see that $H$ is hermitian and essentially self-adjoint on $\Delta$, and also that it commutes with the inversion $I$. So $B = IAI + H = -A + H$ and we have 

{\bf Theorem II ([Bu98a,b]):}\quad The dilation invariant ``conductor operator'' $A + B$ (with initial domain $\Delta$) is essentially self-adjoint, bounded below, commutes with the inversion. Its spectral function is given by the logarithmic derivative of the Tate-Gel'fand-Graev Gamma function on the critical line $\;\bullet$

It is clear that ${\cal S}_0(K)$ is a core for $H$, but the larger domain $\Delta$ is more useful as it is stable under $A, B, I, \cal F$.

{\bf The scattering matrix}

In the Lax-Phillips scattering theory [LP90] one first has a separable Hilbert space. We will take it to be $L^2(K_\nu,dx)$. Then one needs a unitary representation of the additive real line, we will replace that by the unitary representation $U(t)$ of $K_\nu^\times$ through dilations on $L^2$. The next data will be a closed subspace ${\cal D}_-$ of ``Cauchy data for incoming waves'' such that
$$|t|\leq 1 \Rightarrow U(t){\cal D}_-\subset{\cal D}_-$$
$$\bigwedge U(t){\cal D}_- = 0\quad\overline{\bigvee U(t){\cal D}_-} = L^2$$
and a subspace ${\cal D}_+$ of ``Cauchy data for outgoing waves'' such that
$$|t|\geq 1 \Rightarrow U(t){\cal D}_+\subset{\cal D}_+$$
$$\bigwedge U(t){\cal D}_+ = 0\quad\overline{\bigvee U(t){\cal D}_+} = L^2$$
For the theory to gain additional momentum one needs one more axiom
$${\cal D}_-\perp{\cal D}_+$$
The standard example is then of course to take ${\cal D}_-$ to be
$${\cal D}_-^0 = \left\{\varphi(x) \in L^2: |x| > 1 \Rightarrow \varphi(x) =0\right\}$$
and ${\cal D}_+$ to be 
$${\cal D}_+^0 = \left\{\varphi(x) \in L^2: |x| \leq 1 \Rightarrow \varphi(x) =0\right\}$$

{\bf Lemma:} Assume that the unitary operator $\alpha$ of $L^2$ commutes with the $U(t)$'s and restricts to an isometry of ${\cal D}_-^0$ with itself. Then $\alpha$ is a multiple of the identity in each isotypical component $L^2_\chi$ $\;\bullet$

{\bf Proof:}\quad A well-known result of the theory of Hardy spaces. Looking at the $\chi$ component of the spectral representation one sees that $\alpha$ becomes multiplication by a function of modulus $1$ which has a holomorphic continuation to the ``exterior'' of the unit circle, including the point at infinity. It can not vanish anywhere in this domain (else the inclusion $\alpha({\cal D}_-^0) \subset {\cal D}_-^0$ would be strict), so admits a logarithm. Its real part is harmonic and vanishes on the boundary, hence $\alpha(z)$ is a constant (of modulus $1$) $\;\bullet$

So if we have another incoming space ${\cal D}_-$ and an isometry $\alpha$ from $L^2$ to itself (in the commutant of the dilations) such that $\alpha({\cal D}_-) = {\cal D}_-^0$ then the collection of functions $\alpha(\chi;z)$ is unique up to a multiplicative constant in each component.

Assuming in the same manner that we can find an isometry $\beta$ from $L^2$ to itself identifying a given outgoing space ${\cal D}_+$ with ${\cal D}_+^0$ (and commuting with dilations), the operator $S = \beta\ \alpha^{-1}$ is called the scattering matrix. It is a unitary operator, commuting with the dilations, hence defined by a collection of functions of unit modulus $S(\chi;z)$ defined up to multiplicative constants.

Let us try now to find a scattering matrix involving the Fourier Transform. For this we try first ${\cal D}_- = {\cal D}_-^0$ and ${\cal D}_+ = {\cal F(D_-)}$. The scattering matrix is then $U({1\over \pi})({\cal F}I)^{-1}$. But the last axiom is not satisfied because of the characteristic function $\omega$ of the integers which is not perpendicular to its Fourier Transform. Fortunately this is the only obstruction:

{\bf Theorem III:}\quad Let ${\cal D}_- = \left\{\varphi \in {\cal D}_-^0: \int \varphi(x)dx = 0\right\}$. Then ${\cal D}_- \perp {\cal F(D_-)} \; \bullet$

{\bf Proof:}\quad The Fourier Transform of an element $\varphi$ of ${\cal D}_-^0$ is constant on balls of radius $q^\delta$ (as $\lambda(xy) = 1$ for $|x|\leq q^\delta$ and $|y|\leq 1$). If $\varphi$ belongs to ${\cal D}_-$, its Fourier transform vanishes at the origin and is thus supported outside of the ball of radius $q^\delta$ centered at the origin $\;\bullet$

In the remaining sections of this paper ${\cal D}_-$ will always refer to the space just defined and ${\cal D}_+$ to $\cal F(D_-)$. To compute the scattering matrix we first need to find a unitary $\alpha$ commuting with dilations and such that $\alpha({\cal D}_-) = {\cal D}_-^0$. Let us work directly in the spectral picture. There is nothing to do on the components corresponding to ramified characters, we can take $\alpha(\chi;z) = 1$ there. On the component containing the trivial character, a little computation shows that elements of ${\cal D}_-$ are those elements of ${\cal D}_-^0$ for which the associated holomorphic function in the exterior domain vanishes at $z = \sqrt{q}$. So we can take
$$\alpha(1;z) = {\sqrt{q}z - 1\over \sqrt{q} - z}$$
The operator $\beta$ can then be taken as $U({1\over \pi})I\alpha{\cal F}^{-1}$. As $I\alpha I = \alpha^{-1}$ and as ${\cal F}I = \Gamma$ commutes with $\alpha$ this gives as scattering matrix $S = U({1\over \pi})\alpha^{-2}\Gamma^{-1}$. Using the tabulated values of $\Gamma(\chi;z)$ (Theorem I) and the freedom to adjust by multiplicative constants to get rid of the Gauss sums, we end up with:

{\bf Theorem IV:}\quad The scattering matrix $S(\chi;z)$ is given as
$$\hbox{for $\chi$ trivial:}\;\ S(1;z) = z^{-\delta}\ {{\sqrt{q}-z}\over{\sqrt{q}z-1}}$$
$$\hbox{for $\chi$ ramified:}\;\ S(\chi;z) = z^{1 - e(\chi) - \delta}$$

We note that the scattering matrix is ``causal'' in the sense that its spectral function is holomorphic in the exterior domain, having (for $\chi = 1$) a simple zero at $\sqrt{q}$ and a zero of multiplicity $\delta$ at $\infty$. It also has a simple pole at $1\over \sqrt{q}$ and pole of multiplicity $\delta$ at the origin (possibly $\delta = 0$).

{\bf The interacting space and its associated kernel}

The interacting space $\cal K$ is defined as $L^2 - (\cal D_- \oplus D_+)$. We denote its associated orthogonal projector by $K$. The isotypic component $\cal K_\chi$ for a ramified character is easily determined, using the formulae for the $\Gamma$ functions: it has dimension $e(\chi) + \delta - 1$ and is spanned by the functions on $X_\chi$ $e_j(\chi;z) = z^j, 1\leq j\leq e(\chi) + \delta - 1$ . The component ${\cal K}_1$ corresponding to the trivial character contains at least (according to the proof of Theorem III) the function $\omega$ as well as $e_j(1;z)\ 1\leq j\leq \delta $. On the other hand the operator $\beta = S\alpha$ which identifies ${\cal D}_+$ with ${\cal D}_+^0$ is just multiplication by $z^{-\delta}$ on $X_1$, which means that a basis of ${\cal D}_{+,1}$ is given by the functions $e_j(1;z)$ for $j>\delta$. So we have the full count. Let's add the Kronecker delta as a new notation $\delta(\chi_1, \chi_2) = 1$ if $\chi_1 = \chi_2$, $0$ if not.

{\bf Theorem V:}\quad $\cal K = \bigoplus_\chi K_\chi$ with dim($\cal K_\chi$) = $\delta + e(\chi) -1 + 2\delta(\chi, 1) \;\bullet$

We now turn to the representation of the projector $K$ as an integral kernel. For this we just need to find the representation of $q^{\delta/4}\omega(x)$ in the spectral picture. First going from the additive to the multiplicative picture one replaces $\omega(x)$ by $a\sqrt{|t|}\omega(t)$. So $\widehat{\omega}(z) = a\left(1 - {\overline{z}\over\sqrt{q}}\right)^{-1}$. Recalling that $a^2 = q^{-\delta/2}(1 - 1/q)$, we thus deduce:

{\bf Theorem VI:} 
$$K((\chi_1;z), (\chi_2;w)) = \delta(\chi_1, 1)\delta(\chi_2, 1)\left({(1 - 1/q)\over(1 - \overline{z}/\sqrt{q})(1 - w/\sqrt{q})}+ 1_{\delta\geq 1}z^\delta\overline{w}^\delta\right)$$
$$+ \sum_{\chi} \delta(\chi_1, \chi)\delta(\chi_2, \chi)
\sum_{1\leq j \leq \delta + e(\chi) -1} z^j\overline{w}^j$$
$$K((\chi;z), (\chi;z)) = \delta + \delta(\chi, 1){(1 - 1/q)\over|1 - z/\sqrt{q}|^2}+ \sum_{\chi_1\neq 1} \delta(\chi, \chi_1) (e(\chi) -1)\;\bullet$$

{\bf The time delay operator and the conductor operator}

First a little lemma from the the theory of Fourier Series:

{\bf Well-known Kernel Lemma:}\quad Let $A(z,w)$ be a smooth kernel on $S^1 \times S^1$.  Then the series $$\sum_{j \in \ZZ} \int\int z^{-j}A(z,w)w^j\;{d\theta(z)\over 2\pi}{d\theta(w)\over 2\pi}$$ converges absolutely and has value $\int A(z,z) {d\theta(z)\over 2\pi}$ $\;\bullet$

{\bf Proof:}\quad The $\int\int z^{-j}A(z,w)w^j\;{d\theta(z)\over 2\pi}{d\theta(w)\over 2\pi}$ are the Fourier coefficients of the function $\phi(z) = \int A(zw,w)\,{d\theta(w)\over 2\pi}$. Alternatively this is the well-known evaluation of the Hilbert Trace of the operator with integral kernel $A(z,w)\;\bullet$

Let's abbreviate the dilation $U({1\over \pi})$ by just $U$. In the spectral represntation $U$ acts as multiplication by $z$ in each component. The projector $K$ commutes with the units. Lax-Phillips [LP78] define the time delay operator $T$ to measure the average amount of time spent in the interaction with the scatterer. In our multiplicative picture it is natural for the measure of this time to assign a weight of $\log(q)$ to the units. As $<Kf|f>$ is the probability of interaction at ``time in the units'', the total expected time spent interacting with the scatterer is thus $\log(q) \sum_{j \in \ZZ} <KU^jf|U^jf>$, if this makes sense. According to the Lemma this does make sense if $f$ is taken from the domain $\Delta$ considered before, and its value is then obtained as $<Tf|f>$ where $T$ is the dilation invariant operator with spectral function $T(\chi;z) = \log(q) K((\chi;z), (\chi;z))$. Obviously $(T, \Delta)$ is essentially self-adjoint and non-negative. Let us record this as a theorem (its last statement is checked with an easy explicit computation):

{\bf Theorem VII:}\quad The dilation invariant, non-negative, essentially self-adjoint operator $T$ on $\Delta$ with spectral function
$$T(\chi;z) = \log(q)\left(\delta + \delta(\chi, 1){(1 - 1/q)\over|1 - z/\sqrt{q}|^2}+ \sum_{\chi_1\neq 1} \delta(\chi, \chi_1) (e(\chi) -1)\right)$$
is such that 
$$\{f \in \Delta\} \Rightarrow <Tf|f> = \int_{K_\nu^\times} ||KU(t)f||^2 \,d^\times t$$
and is related to the scattering matrix through
$$T(\chi;z) = S(\chi;z)\, D(\overline{S(\chi;z)})\quad\bullet$$

As we have both the time delay operator and the conductor operator explicitely we can relate them easily:

{\bf Theorem VIII:}\quad $T = \delta\log(q) + |(\delta + 1) \log(q) - H|\quad\bullet$

{\bf Proof:}\quad one just has to compare the spectral function of $T$ given above with the spectral function of $H$ given before Theorem II: one finds
$$\hbox{on }L^2_1\quad0\leq T - \delta\log(q) = - H + (\delta + 1)\log(q)$$
$$\hbox{on }L^2_\chi\hbox{ for }\chi\hbox{ ramified}\quad 0\leq T - \delta\log(q) = H - (\delta + 1)\log(q)\;\bullet$$

{\bf The contraction semi-group and a trace formula}

Let $P_-$ be the orthogonal projection onto ${\cal K} \oplus {\cal D}_+$, and $P_+$ be the orthogonal projection onto ${\cal K} \oplus {\cal D}_-$. As $U(t)$ leaves stable ${\cal D}_-$ for $|t|\leq1$ the vectors $U(t)f$, $|t|\geq1$ remain orthogonal to ${\cal D}_-$ for any $f$ in ${\cal K} \oplus {\cal D}_+$. On the other hand the subspace ${\cal D}_+$ is stable under $U(t)$, $|t|\geq1$. So we get a semi-group of operators $Z(t)$, $|t|\geq1$ on the quotient ${\cal K}$. The operators are given by the formula $Z(t) = P_+U(t)P_-$, for $|t|\geq1$.

For $t$ a unit $Z(t)$ is just multiplication by $\chi(t)$ in the isotypical component ${\cal K}_\chi$. $Z({1\over\pi})$ is multiplication by the function ``$z$'' on each isotypical component ${\cal K}_\chi$. It is thus nilpotent on the orthogonal space to Tate's function $\omega$ and $<Z({1\over\pi})\omega|\omega> = <U({1\over\pi})\omega|\omega> = {1\over\sqrt{q}}<\omega|\omega>$ so that its eigenvalue is ${1\over\sqrt{q}}$.

Let us define $Z(f)$ for $f(t)$ a Schwartz-Bruhat function on $K_\nu^\times$ (with support in $|t|\geq 1$) as
$$Z(f) = \int_{|t|\geq 1} f(t)Z(t) \,d^\times t$$

The operator $Z(f)$ is of trace-class, in fact it is of finite rank: one first writes $f$ as a finite linear sum $f = \sum_\chi f_\chi$ with $f_\chi = \sum_{j\geq0} U({1\over\pi})^j(f_{\chi,j})$, where each $f_{\chi,j}$ has support on the units and is a multiple of $\chi(1/t)$ there. As $Z(t) = \chi(t)$ on ${\cal K}_\chi$ for $t$ a unit, $Z(f_\chi)$ acts non-trivially only on ${\cal K}_\chi$ (which is finite dimensional). Using the fact that $U({1\over\pi})$ stabilizes the orthogonal complement to $\omega$ and is nilpotent there, while its eigenvalue on the $1-$dimensional quotient is ${1\over\sqrt{q}}$, one evaluates the traces as follows
$$\hbox{Tr}(Z(f_1)) = \left(\int_{|t|=1} f_1(t)\,d^\times t\right)\hbox{dim}({\cal K}_1) + \int_{|t|>1} |t|^{-1/2}f_1(t)\,d^\times t$$
$$(\chi \neq 1)\;\hbox{Tr}(Z(f_\chi)) = \left(\int_{|t|=1} f_\chi(t)\chi(t)\,d^\times t\right)\hbox{dim}({\cal K}_\chi) $$

{\bf Theorem IX:}\quad Tr$(Z(f)) = T(f)(1) = \int_{X_\nu} \widehat{f}(\chi;z)T(\chi;z)\,{d\theta\over 2\pi}\;\bullet$

{\bf Proof:}\quad First we look at $f_\chi$. The spectral function of $T$ is the constant $(\delta + e(\chi) - 1)\log(q)$ on $X_\chi$ (Theorem VII) while the dimension of ${\cal K}_\chi$ is $(\delta + e(\chi) - 1)$ (Theorem V). Then we consider $f_1$. From Theorem V we know that dim$({\cal K}_1)$ is $\delta + 1$, so we can rewrite Tr$(Z(f_1))$ as $f_1(1)\delta \log(q) + \int_{|t|\geq1}|t|^{-1/2}f_1(t)\,d^\times t$. In terms of the spectral transform $\widehat{f_1}(z)$ which is the boundary value of a holomorphic function in the inner domain of the circle $X_1$, this becomes
$$\left(\int_{X_1} \widehat{f_1}(z)\,{d\theta\over 2\pi}\right)\delta\log(q) + \widehat{f_1}({1\over\sqrt{q}})\log(q)$$
On the other hand the spectral function of $T$ on $X_1$ was given in Theorem VII as 
$$T(1;z) = \left(\delta + {(1 - 1/q)\over|1 - z/\sqrt{q}|^2}\right)\log(q)$$
which is the boundary value of the meromorphic function
$$\left(\delta + {(1 - 1/q)z\over(1 - z/\sqrt{q})(z - 1/\sqrt{q})}\right)\log(q)$$
Integrating against ${d\theta\over 2\pi}$ means a complex line integral against ${1\over2\pi\,i}{dz\over z}$. So applying Cauchy residue formulae gives the stated answer$\;\bullet$ 

{\bf Weil's local term in terms of a supertrace}

We now extend the definition of $Z(f)$ to the case that $f$ has its support in $|t| < 1$ as $Z(f) = Z(If)$, and then by linearity to the general case of a Schwartz-Bruhat function on $K_\nu^\times$. As the time delay operator $T$ commutes with the inversion the trace formula of Theorem IX remains valid.

Using the relation between the operators $T$ and $H$ we will now give a trace formula for $H(f)(1)$. For this purpose let us introduce a ``superspace'' structure on ${\cal K}$ where the odd part will be the one-dimensional space spanned by $\omega$ and the even part its orthogonal complement. The super-dimensions of the isotypical components are
$$\hbox{sDim}({\cal K}_\chi) = \delta + e(\chi) -1$$
We first evaluate sTr$(Z(f))$ for $f$ with support in $\{|t|\geq 1\}$. We just need to go back to the considerations before Theorem IX and one gets
$$\hbox{sTr}(Z(f_1)) = \left(\int_{|t|=1} f_1(t)\,d^\times t\right)
\hbox{sDim}({\cal K}_1) - \int_{|t|>1} |t|^{-1/2}f_1(t)\,d^\times t$$
$$(\chi \neq 1)\;\hbox{sTr}(Z(f_\chi)) = \left(\int_{|t|=1} f_\chi(t)\chi(t)\,d^\times t\right)\hbox{sDim}({\cal K}_\chi)$$
$$\hbox{sTr}(Z(f_1)) + f_1(1)\log(q) = f_1(1)\delta\log(q) - \int_{|t|>1} |t|^{-1/2}f_1(t)\,d^\times t$$
$$(\chi \neq 1)\;\hbox{sTr}(Z(f_\chi)) + f_\chi(1)\log(q) = f_\chi(1)(\delta + e(\chi))\log(q)$$

{\bf Theorem X:}\quad sTr$(Z(f)) + f(1)\log(q) = H(f)(1) = \int_{X_\nu} \widehat{f}(\chi;z)H(\chi;z)\,{d\theta\over 2\pi}\;\bullet$

{\bf Proof:}\quad The spectral function of $H$ is the constant $(\delta + e(\chi))\log(q)$ on $X_\chi$ (Theorem II) so this settles it for $f_\chi$. For the invariant component $f_1$, $$\hbox{sTr}(Z(f_1)) + f_1(1)\log(q) + \hbox{Tr}(Z(f_1)) = f_1(1)(2\delta + 1)\log(q)$$
On the other hand we also have
$$H(1;z) + T(1;z) = (2\delta + 1)\log(q)$$
so the result is then implied by Thorem IX. Finally $H$ commutes with the inversion so the result also holds if $f$ has its support in $\{|t|< 1\}$. The general case then follows by linearity$\;\bullet$ 

According to [Bu98a] the quantity $H(f)(1)$ is the local contribution to the Explicit Formula [We52] (with the incorporation of the local component of the discriminant, and with the convention that the Explicit Formula counts the zeros with positive multiplicity and the poles with negative multiplicity, and with a shift of ${1\over2}$ so that the poles are at $\pm{1\over2}$). Actually the explicit evaluations stated above just before Theorem X are another derivation of this fact.

{\bf Exact evaluation of a trace considered by Connes}

As a further application of the method let's evaluate the simplest of the traces considered by Connes (Theorem $3$ of Section V from [Co98]). Connes considers the operator
$$\widetilde{P_\Lambda}P_\Lambda U(f)$$
where $P_\Lambda$ is orthogonal projection to functions with support in $|x|\leq\Lambda$, $\widetilde{P_\Lambda}$ is its Fourier conjugate and 
$$U(f) = \int_{K_\nu^\times} f(t)U(t) \,d^\times t$$
which is $\log(q)$ times the multiplicative convolution with the locally constant compactly supported function $f$ (in fact Connes's $U(h)$ is our $U(f)$ for $f(t) = {\sqrt{|t|}} h(t)$). He shows (also at an archimedean place) that it is of trace class, has a main logarithmic divergency and a constant term which is the local term of the Explicit Formula and a $o(1)$ error term. In the non-archimedean case his proof actually gives the exact value for $\Lambda$ large enough. 

Let us evaluate this trace with our methods. Put $\Lambda = q^n, n \geq 0$, $P_n = P_\Lambda$, and $$Q_n = \widetilde{P_n}P_n$$
As $P_n$ and $\widetilde{P_n}$ commute with the action of the unit group, $Q_n$ respects the decomposition of the Hilbert space in isotypical components. Let $Q_n^\chi$ be the corresponding operator on $L_\chi^2$. 

{\bf Theorem XI:}\quadÊThe operator $Q_n^\chi$ is, in the spectral picture, an integral operator with a smooth kernel $Q_n^\chi(z,w)$. Its value is as follows for $\chi$ ramified:
$$Q_n^\chi(z,w) = \sum_{e(\chi)+\delta - n\leq j \leq n} z^j\overline{w}^j \quad\hbox{(this vanishes if $2n < e(\chi) + \delta$)}$$
and, for the trivial character:
$$(2n\geq \delta)\quad\quad Q_n^1(z,w) = {(z\overline{w})}^{\delta -n}\left({(1 - 1/q)\over(1 - \overline{z}/\sqrt{q})(1 - w/\sqrt{q})} + \sum_{1\leq j \leq 2n - \delta} {(z\overline{w})}^j\right)$$
$$(0\leq 2n < \delta)\quad\quad Q_n^1(z,w) = {z^{\,\delta - n}\;\overline{w}^{\,n}\over { \sqrt{q}^{\,\delta - 2n}}}{(1 - 1/q)\over(1 - \overline{z}/\sqrt{q})(1 - w/\sqrt{q})}$$
For $n$ given, only finitely many of the kernels $Q_n^\chi(z,w)$ are non-zero.$\;\bullet$

{\bf Proof:}\quad\quad First we consider a ramified character. The action of $P_n$ is orthogonal projection to the span of $\{z^j, j\leq n\}$ while the action of $\widetilde{P_n}$ is by Theorem I orthogonal projection to the span of $\{z^j, j\geq e(\chi) + \delta -n\}$. Hence the result.

For the trivial character, let us first assume that $n\geq \delta - n$. We work exclusively with functions invariant under the units. Let, for $j \in \ZZ$, $\omega_j$ be the function of norm $1$ constant on the ball of radius $q^j$. The Fourier Transform of $\omega_j$ is $\omega_{\delta - j}$. Let ${\cal K}_n$ be the space spanned by $\omega_{\delta - n}, \dots , \omega_{n}$, let ${\cal L}_n$ the space spanned by functions perpendicular to $\omega_{\delta - n}$, with support in $|x| \leq q^{\delta - n}$, and ${\cal M}_n$ the space of functions with support in $|x| > q^{n}$. Then $L^2_1 = {\cal L}_n \oplus {\cal K}_n \oplus {\cal M}_n$. $P_n$ cuts off ${\cal M}_n$, $\widetilde{P_n}$ cuts off ${\cal L}_n$, so they commute and their combination $Q_n$ is just orthogonal projection to ${\cal K}_n$. An orthonormal basis of this space in the spectral picture is given by the functions $\sqrt{1 - {1\over q}}z^{\delta - n} \left(1 - {\overline{z}\over\sqrt{q}}\right)^{-1}$ and $\{ z^j, \delta - n < j \leq n\}$. So
$$Q_n^1(z,w) = {(z\overline{w})}^{\delta - n}\left({(1 - 1/q)\over(1 - \overline{z}/\sqrt{q})(1 - w/\sqrt{q})} + \sum_{1\leq j \leq 2n - \delta} {(z\overline{w})}^j\right)$$
If $n < \delta - n$, let $m = \delta - n$ so that $m\geq \delta - m$. Then $P_n$ is identically $1$ on ${\cal L}_m$ and identically zero on ${\cal M}_m$, and $\widetilde{P_n}$ is identically $1$ on ${\cal M}_m$ and identically zero on ${\cal L}_m$, so the combination $Q_n$ kills both ${\cal L}_m$ and ${\cal M}_m$. On ${\cal K}_m$, $P_n$ sends a given $\phi$ to $<\omega_n | \phi> \omega_n$ which is then sent by $\widetilde{P_n}$ to $<\omega_n | \phi> <\omega_n | \omega_m>\omega_m = \left({1\over \sqrt{q}}\right)^{\delta - 2n}<\omega_n | \phi> \omega_m$. This gives the kernel
$$Q_n^1(z,w) = {z^{\,\delta - n}\;\overline{w}^{\,n}\over { \sqrt{q}^{\,\delta - 2n}}}{(1 - 1/q)\over(1 - \overline{z}/\sqrt{q})(1 - w/\sqrt{q})}\eqno{\bullet}$$

{\bf Theorem XII:}\quad For any $n\geq 0$ and any $f$ in the domain $\Delta$ the operator $Q_n U(f)$ is a finite rank operator, whose Hilbert Trace is
$$\hbox{Tr}(Q_n U(f)) = \log(q) \sum_\chi \int_{X_\chi} Q_n^\chi(z,z) \widehat{f}(\chi;z) {d\theta(z) \over 2\pi}$$ 

{\bf Proof:}\quad If $f \in \Delta$ then $\widehat{f}(\chi;z)$ vanishes for almost all $\chi$'s. Of course $U(f)$ acts through multiplication with $\log(q)\widehat{f}(\chi;z)$ in the spectral picture, and as the previous proof shows that each $Q_n^\chi$ has finite rank, $Q_n U(f)$ also has finite rank. Its trace is given by the ``Well-known Trace Lemma'' we proved earlier$\;\bullet$

For $m\geq 0$ let $\epsilon_m(t)$ be the function on $K_\nu^\times$ whose spectral function is the characteristic function of the components indexed by characters with conductors at most $m$. For $f \in \Delta$  the multiplicative convolution $\epsilon_m * f$ coincides with $f$ for $m \geq e(f)=$ the largest conductor from the $\chi$'s for which $\widehat{f}(\chi;z)$ is non-zero.

{\bf Theorem XIII:} \quad For $f \in \Delta$:
$$\hbox{Tr}(Q_n U(f)) = (2n + 1) (\epsilon_{2n - \delta} * f)(1) \log(q) - H(\epsilon_{2n - \delta} * f)(1)\leqno{(2n \geq \delta)}$$
$$\hbox{Tr}(Q_n U(f)) = (2n + 1) f(1) \log(q) - H(f)(1)\leqno{(2n \geq \delta + e(f))}$$

{\bf Proof:}\quadÊOne just need to compare $Q_n^\chi(z,z)$ to the spectral function of the conductor operator $H(\chi; z)$. One finds $Q_n^\chi(z,z)\log(q) = (2n + 1)\log(q) - H(\chi; z)$ for $2n \geq e(\chi) + \delta$ and $Q_n^\chi(z,z) = 0$ for $\delta\leq 2n < \delta + e(\chi)$ (which happens only for a ramified character)$\quad\bullet$

{{\bf REFERENCES}\par
\baselineskip = 12 pt
\parskip = 4 pt
\font\smallRoman = cmr8
\smallRoman
\font\smallBold = cmbx8
\font\smallSlanted = cmsl8
{\smallBold [Bu98a] J.F. Burnol}, {\smallSlanted ``The Explicit Formula and a Propagator''}, math/9809119 (v1 Sep{. }1998, v2 November 1998)\par
{\smallBold [Bu98b] J.F. Burnol}, {\smallSlanted ``Spectral analysis of the local conductor operator''}, math/9811040 (November 1998)\par
{\smallBold [Bu98c] J.F. Burnol}, {\smallSlanted ``Spectral analysis of the local commutator operators''}, math/9812012 (December 1998)\par
{\smallBold [Co98] A. Connes}, {\smallSlanted ``Trace formula in non-commutative Geometry and the zeros of the Riemann zeta function''}, math/9811068 (November 1998).\par
{\smallBold [LP78] P.D. Lax, R.S. Phillips}, {\smallSlanted ``The time delay operator and a related trace formula''} in {\smallSlanted ``Topics in Functional Analysis''} (I. Gohberg and M.Kac, eds.), pp 197-215, Academic Press, New York (1978).\par
{\smallBold [LP90] P.D. Lax, R.S. Phillips}, {\smallSlanted ``Scattering Theory''}, Academic Press, San Diego (1990, 2nd ed.).\par
{\smallBold [Ta50] J. Tate}, {\smallSlanted ``Fourier Analysis in Number Fields and Hecke's Zeta Function''}, Princeton 1950, reprinted in Algebraic Number Theory, ed. J.W.S. Cassels and A. Fr\"ohlich, Academic Press, (1967).\par
{\smallBold [We52] A. Weil},{\smallSlanted ``Sur les ``formules explicites'' de la th\'eorie des nombres premiers''}, Comm. Lund (vol d\'edi\'e \`a Marcel Riesz), (1952).

\vfill
\centerline{Jean-Fran\c{c}ois Burnol}
\centerline{62 rue Albert Joly}
\centerline{F-78000 Versailles}
\centerline{France}
\centerline{January 1999}
}
\eject
\bye